\documentclass[11pt]{article}
\usepackage{amssymb}
\usepackage{amsmath}
\usepackage{pst-all}
\usepackage{graphicx}
\usepackage{pstricks}
\usepackage{makeidx}
\usepackage{enumerate}
\usepackage{mathtools}
\usepackage{comment}
\usepackage{hyperref}
\usepackage{tabularx}
\usepackage[a4paper,left=2cm,right=2cm,top=2.5cm,bottom=2.5cm]{geometry}

\usepackage{tikz}
\usetikzlibrary{arrows.meta, positioning}
\usetikzlibrary{decorations.pathreplacing}
\usetikzlibrary{calc,positioning}

\usetikzlibrary{trees,shapes}

\usepackage{tabularx} 
\usepackage{booktabs}

\usepackage{multirow}
\usepackage{graphicx}

\newtheorem{lemma}{Lemma}[section]
\newtheorem{theorem}[lemma]{Theorem}
\newtheorem{corollary}[lemma]{Corollary}
\newtheorem{question}[lemma]{Question}
\newtheorem{proposition}[lemma]{Proposition}
\newtheorem{remark}[lemma]{Remark}
\newtheorem{definition}[lemma]{Definition}
\newtheorem{example}[lemma]{Example}
\newtheorem{problem}[lemma]{Problem}

\def\bnum{\begin{enumerate} }
\def\enum{\end{enumerate}}
\def\bdf{\begin{definition}\rm }
\def\edf{\end{definition}}
\def\br{\begin{remark}\rm }
\def\er{\end{remark}}
\def\be{\begin{equation}}
\def\ee{\end{equation}}
\def\bt{\begin{theorem}}
\def\et{\end{theorem}}
\def\bl{\begin{lemma}}
\def\el{\end{lemma}}
\def\bc{\begin{corollary}}
\def\ec{\end{corollary}}
\def\bp{\begin{proposition}}
\def\ep{\end{proposition}}
\def\bxa{\begin{example}\rm }
\def\exa{\end{example}}
\def\ba{\begin{array}}
\def\ea{\end{array}}
\def\ben{\begin{eqnarray*}}
\def\een{\end{eqnarray*}}
\def\bdsc{\begin{description}}
\def\edsc{\end{description}}
\def\bpsp{\begin{pspicture}}
\def\epsp{\end{pspicture}}
\def\bea{\begin{eqnarray}}
\def\eea{\end{eqnarray}}
\def\btab{\begin{tabular}}
\def\etab{\end{tabular}}
\def\bpm{\begin{problem}}
\def\epm{\end{problem}}
\def\bfig{\begin{figure}}
\def\efig{\end{figure}}

\def\bnum{\begin{enumerate}\itemsep=0cm}
\def\enum{\end{enumerate}}

\def\pr{{\em Proof. }}

\def \qe{\hfill \vrule height4pt width 4pt depth 0pt}

\def\1{1\!\hspace{-.08cm}1}

\title{On the Seidel energy of uniform hypergraphs due to hyperedge and vertex deletion} 

\author{Shib~Sankar~Saha\thanks{Corresponding author: Department of Mathematics, Bar-Ilan University, Israel (sahashi@biu.ac.il)}
}
\begin{document}
\pagestyle{myheadings}
\markboth{S.~S.~Saha}{\title{On the Seidel energy of uniform hypergraphs due to hyperedge and vertex deletion}}
\maketitle
\begin{abstract}
Let $\mathcal{S}(\mathcal{H})$ be the Seidel matrix of a hypergraph $\mathcal{H}$, and the Seidel energy is denoted by the sum of the absolute eigenvalues of $\mathcal{S}(\mathcal{H})$. In [G.~X.~Tian, Y.~Li and  S.~Y.~Cui, The change of Seidel energy of tripartite Turán graph due to edge deletion, Linear Multilinear Algebra, 19 (2022), 4597-4614] and [Y.~Liu, X.~Chen, The change of Seidel energy of 5-partite Turán graph due to edge deletion, Discrete Applied Mathematics, 2024, 342, 104-123], the authors studied the change of Seidel energy of the Turán graph due to edge deletion. In this article, we analyze the Seidel spectrum of the complete $3$-uniform bipartite hypergraph $\mathcal{C}^3_{m,n}$ and show that it has exactly one negative Seidel eigenvalue even after a single hyperedge deletion. Finally, we prove that the Seidel energy of the complete $3$-uniform bipartite hypergraph $\mathcal{C}^3_{m,n}$ decreases after single hyperedge and vertex deletion for all $m,n \ge 3$.
\end{abstract}

\noindent{\bf Keywords.}  Hypergraph; Seidel matrix; Equitable partition; Seidel spectrum; Seidel energy.

\noindent{\bf Mathematics Subject Classifications.} 05C65; 05C35; 05C50; 15A18; 15A60
\sloppy

\section{Introduction}
Throughout this paper, all hypergraphs are assumed to be finite, simple and connected. A hypergraph $\mathcal{H}=\big(V(\mathcal{H}), E(\mathcal{H})\big)$ consists of a vertex set $V(\mathcal{H})$ and a hyperedge set $E(\mathcal{H})\subseteq 2^{V(\mathcal{H})}$. For any $k\in\mathbb{Z}^+$ and any set $X$, let $[k]=\{1,2,\ldots,k\}$, $\binom{X}{k}=\{Y\subseteq X: |Y|=k\}$ and $\binom{|X|}{k}$ is called a binomial coefficient and it represents the number of ways to choose $k$ elements from a set X of size $|X|$, without regard to order. For a positive integer $k$, a hypergraph $\mathcal{H}$ is said to be $k$-uniform if $E(\mathcal{H})\subseteq\binom{V(\mathcal{H})}{k}$, and a $k$-uniform hypergraph is also called a $k$-graph. Especially, a graph can be regarded as a $2-$uniform hypergraph. For an undirected hypergraph, each hyperedge is an unordered set. A hypergraph $\mathcal{H}$ is {\em finite} if $V(\mathcal{H})$ is a finite set. A hypergraph $\mathcal{H}$ is said to be simple if it does not contain any hyperedge e such that $|e|=1$ and for any $e_i,e_j\in E(\mathcal{H})$ such that $e_i\neq e_j$ and neither $e_i\subset e_j$ nor $e_j\subset e_i$. Let $\mathcal{H}'= \big(V(\mathcal{H}'), E(\mathcal{H}')\big)$ be a hypergraph with $V(\mathcal{H}')\subseteq V(\mathcal{H})$ and $E(\mathcal{H}')\subseteq E(\mathcal{H})$, then $\mathcal{H'}$ is a subhypergraph of $\mathcal{H}$.

The \emph{complete $k$-uniform bipartite hypergraph} $\mathcal{C}^{k}_{t,s}=\big(V(\mathcal{C}^{k}_{t,s}), E(\mathcal{C}^{k}_{t,s})\big)$ is the $k$-uniform hypergraph on a bipartition $V(\mathcal{C}^{k}_{t,s})=V_1\sqcup V_2$ with $|V_1|=m$, $|V_2|=n$, whose hyperedges are
\[
E(\mathcal{C}^{k}_{t,s})
 = \{\, e\subseteq V_1\cup V_2: |e|=k,\ e\cap V_1\neq\varnothing,\ e\cap V_2\neq\varnothing \,\}.
\]

The notation $\mathcal{H}-e=\big(V(\mathcal{H}-e), E(\mathcal{H}-e)\big)$ represents a {\em deletion (weak deletion)} of a hyperedge $e\in E(\mathcal{H})$ from $\mathcal{H}$ such that $V(\mathcal{H}-e)=V(\mathcal{H})$, $E(\mathcal{H}-e)=E(\mathcal{H})\setminus\{e\}$ and is the simplest operation of deletion in a hypergraph. We use $i\sim j$ to denote that the vertex $i$ is adjacent to the vertex $j$. A hypergraph $\mathcal{H}$ is said to be {\em connected} if for any two distinct vertices $i,j$ there exists a sequence of hyperedges $e_{1},e_{2},\ldots,e_{m}$ such that $i\in e_{1}$, $j\in e_{m}$ and $e_{r}\cap e_{r+1}\neq \emptyset$  $\forall~1\leq r\leq m-1$. For a hypergraph $\mathcal{H}$, we use the notation $E_{ij}=\{e\in E(\mathcal{H}): i,j\in e\}$ and $c_{ij}$ is the {\em co-degree} of $\mathcal{H}$, i.e., the cardinality of the set $E_{ij}$ and denoted by $|E_{ij}|$. For more details, see \cite{bretto2013, vol2009}.

Let $A(G)$ be the adjacency matrix of $G$. The Seidel matrix of $G$ is the matrix $S(G) = J_n-I_n-2A(G)$, where $J_n$ denotes the all the entry being $1$ of size $n\times n$ and $I_n$ denotes the identity matrix of order $n$. Thus, the $(i, j)$-entry of $S(G)$ is $0$ if $i = j$, $-1$ if $i$ adjacent to $j$, and $1$ otherwise. The introduction of Seidel matrices is significant as a technique for investigating equiangular line systems in Euclidean spaces \cite{S1996}.

Let $M$ be a Hermitian matrix of order $n$. Then $M$ has $n$ real eigenvalues, which can be arranged in non-increasing order, i.e., $\lambda_1(M),\lambda_2(M),\ldots,\lambda_n(M)$. The energy of $M$, denoted by $\mathcal{E}(M)$, is defined to be the sum of the absolute values of the eigenvalues of $M$. Let $\lambda_1(G),\lambda_2(G),\ldots,\lambda_n(G)$ be the Seidel eigenvalues of $G$. In 2012, W.~H.~Haemers \cite{WH2012} first defined the Seidel energy $\mathcal{E}_S(G)$ of a graph G, which is the sum of absolute values of all Seidel eigenvalues of $G$ i.e., 
\begin{center}
$\mathcal{E}_S(G)=\sum\limits_{i=1}^n|\lambda_i(G)|$,    
\end{center}
and obtained a tight upper bound on the Seidel energy of graphs with n vertices. In 2018, F.~Szöllősi and P.~R.~J.~Östergård \cite{F2018} studied the enumeration of Seidel matrices. In 2019, A.~Berman, N.~Shaked-Monderer, and R.~Singh \cite{A2019} studied complete multipartite graphs that were determined up to switching their Seidel spectrum. In 2019, D.~Rizzolo \cite{D2019} studied the determinants of Seidel matrices and a conjecture of Gurbani. In 2020, S.~Akbari, M.~Einollahzadeh, and M.~M.~Arkhaneei \cite{S2020} studied the proof of a conjecture on the Seidel energy of graphs. In graph theory, the influence of edges on the graph Seidel energy has aroused a large interest among scholars. In 2024, M. Einollahzadeh and M.~M.~Arkhaneei \cite{E2024} proposed a short proof of Haemers' conjecture on the Seidel energy of graphs.

In 1996, K. Feng, W. Ch'ing, and W. Li defined the adjacency matrix $\mathcal{A}({\mathcal{H}})=[\mathcal{A}({\mathcal{H}})_{ij}]_{n\times n}$ of a hypergraph $\mathcal{H}=\big(V(\mathcal{H}), E(\mathcal{H})\big) $ with $n$ vertices as:

$$\mathcal{A}({\mathcal{H}})_{ij}=
\begin{cases}
c_{ij}~~~~,\mbox{ if }~i\sim j,\\
~~0~~~~,\mbox{ otherwise.}
 \end{cases}$$

where $c_{ij}$ is the co-degree of $\mathcal{H}$ \cite{F1996}. 

In this article, we study the concept of the Seidel matrix for a hypergraph. Therefore, the usual definition of the Seidel matrix $\mathcal{S}(\mathcal{H})$ of $\mathcal{H}$ of order $n$ is defined as $\mathcal{S}(\mathcal{H})=J_n-I_n-2\mathcal{A}(\mathcal{H})$, where $\mathcal{A}(\mathcal{H})$, $J_n$, and $I_n$ are the adjacency matrix of a hypergraph $\mathcal{H}$, the matrix with all entries $1$, and the identity matrix, respectively. The entries of the Seidel matrix $\mathcal{S}(\mathcal{H})=[\mathcal{S}(\mathcal{H})_{ij}]_{n\times n}$ are as follows:
$$\mathcal{S}(\mathcal{H})_{ij}=
\begin{cases}
1-2c_{ij}~,\mbox{ if }~i\neq j,\\
~~0~~~~~~~~,\mbox{ if }~i=j.
\end{cases}$$

The Seidel energy of a hypergraph $\mathcal{H}$ is denoted by $\mathcal{E_S}(\mathcal{H})$, which is the sum of absolute values of all Seidel eigenvalues of $\mathcal{H}$, i.e.,
\begin{center}
$\mathcal{E_S}(\mathcal{H})=\sum\limits_{i=1}^n|\rho_i(\mathcal{H})|$,    
\end{center}

where $\rho_1(\mathcal{H}),\rho_2(\mathcal{H}),\cdots,\rho_n(\mathcal{H})$ are the eigenvalues of $\mathcal{S}(\mathcal{H})$ are called the Seidel eigenvalues of $\mathcal{H}$ \cite{LJK2023}. Since $\mathcal{S}(\mathcal{H})$ is a real symmetric matrix, then $\rho_i(\mathcal{H})~(i=1,2,\ldots,n)$ are all real. The multi-set of all eigenvalues of $\mathcal{S}(\mathcal{H})$ is called the Seidel spectrum of $\mathcal{H}$.

\section{Monotonicity of the Seidel energy under hyperedge deletion}
In 2022, G.~X.~Tian, Y.~Li, and S.~Y.~Cui \cite{G2022} studied the change in Seidel energy of the tripartite Turán graph due to edge deletion. In 2024, Y.~Liu, X.~Chen \cite{YL2024} studied the change of Seidel energy of a 5-partite Turán graph due to edge deletion.

Motivated by the above effect of edge deletion on the Seidel energy of graphs, we demonstrate that removing an arbitrary single hyperedge from a hypergraph can affect the Seidel energy of hypergraphs. Deleting an arbitrary single hyperedge may cause the Seidel energy of hypergraphs to increase, decrease, or remain unchanged.

\begin{enumerate}
\item \textbf{Seidel energy increases for an arbitrary single hyperedge deletion:} 
    
Consider the hypergraph $\mathcal{H}=\big(V(\mathcal{H}), E(\mathcal{H})\big)$ such that the vertex set $V(\mathcal{H})=\{1,2,3,4,5,6\}$ and the hyperedge set $E(\mathcal{H})=\Big\{\{1,2,3\},\{1,4,5\},\{2,5,6\},\{3,4,6\}\Big\}$. Then the Seidel spectra of $\mathcal{H}$ and $\mathcal{H}-e$ are $\Big\{-3,-1,-1,-1,3,3\Big\}$ and $\Big\{-\sqrt{5},-\sqrt{5},-\sqrt{5},\sqrt{5},\sqrt{5},\sqrt{5}\Big\} $. For each $e\in E$, we have
$\mathcal{E_S}(\mathcal{H}-e)=13.4164079$ and $\mathcal{E_S}(\mathcal{H})=12$. Hence $\mathcal{E_S}(\mathcal{H}-e)>\mathcal{E_S}(\mathcal{H})$.

\item \textbf{Seidel energy decreases for an arbitrary single hyperedge deletion:}

Consider the hypergraph $\mathcal{H}=\big(V(\mathcal{H}), E(\mathcal{H})\big)$ with the vertex set $V(\mathcal{H})=\{1,2,3,4,5,6,7,8\}$ and the hyperedge set $E(\mathcal{H})=\Big\{\{1,2,3\},\{3,4,5\},\{5,6,7\},\{1,7,8\}\Big\}$. Then the Seidel spectra of $\mathcal{H}$ and $\mathcal{H}-e$ are $\Big\{-1-2\sqrt{2},-1-2\sqrt{2},-1,-1,3,3,-1+2\sqrt{2},-1+2\sqrt{2}\Big\}$ and $\Big\{-1-2\sqrt{3},1-2\sqrt{2},-3,1,1,1,-1+2\sqrt{3},1+2\sqrt{2}\Big\} $. For $e\in E$, we have
$\mathcal{E_S}(\mathcal{H}-e)=18.5850575$ and $\mathcal{E_S}(\mathcal{H})=19.3137085$. Hence $\mathcal{E_S}(\mathcal{H}-e)<\mathcal{E_S}(\mathcal{H})$.

\item  \textbf{Seidel energy equals for an arbitrary hyperedge deletion:}

Consider the hypergraph $\mathcal{H}=\big(V(\mathcal{H}), E(\mathcal{H})\big)$ with the vertex set $V(\mathcal{H})=\{1,2,\ldots,n\}$ and the hyperedge set $E(\mathcal{H})=\Big\{\{1,2,\ldots,n\}\Big\}$. Then the Seidel spectrum of $\mathcal{H}$ are $-(n-1)$ and $1$ with multiplicity $1$ and $n-1$, respectively. Also, the Seidel spectrum of $\mathcal{H}-e$ are $(n-1)$ and $-1$, with multiplicity $1$ and $n-1$, respectively. Then $\mathcal{E_S}(\mathcal{H}-e)=2(n-1)=\mathcal{E_S}(\mathcal{H})$. Hence $\mathcal{E_S}(\mathcal{H}-e)=\mathcal{E_S}(\mathcal{H})$.
\end{enumerate}
These observations naturally motivate the following problems:

\rm\begin{enumerate}
\item [][\textbf{Problem 1}] For which classes of hypergraphs $\mathcal{H}=\big(V(\mathcal{H}), E(\mathcal{H})\big)$ does there exist a hyperedge $e\in E(\mathcal{H})$ whose removal strictly decreases the Seidel energy, that is, $\mathcal{S}(\mathcal{H})>\mathcal{S}(\mathcal{H}-e)?$

\item [][\textbf{Problem 2}] For which classes of hypergraphs$\mathcal{H}=\big(V(\mathcal{H}), E(\mathcal{H})\big)$ does there exist a hyperedge $e\in E(\mathcal{H})$ whose removal leaves the Seidel energy unchanged, that is, $\mathcal{S}(\mathcal{H})=\mathcal{S}(\mathcal{H}-e)$?

\item [][\textbf{Problem 3}] For which classes of hypergraphs $\mathcal{H}=\big(V(\mathcal{H}), E(\mathcal{H})\big)$ does there exist a hyperedge $e\in E(\mathcal{H})$ whose removal strictly increases the Seidel energy, that is, $\mathcal{S}(\mathcal{H})<\mathcal{S}(\mathcal{H}-e)?$
\end{enumerate}

\subsection{Seidel spectrum of $\mathcal{C}^3_{m,n}$ and $\mathcal{C}^3_{m,n}-e$.}
In this section, we address \textbf{Problem 1} by examining the existence of classes of $k$-uniform hypergraphs whose Seidel energy decreases upon the deletion of a hyperedge.

\begin{definition}\label{d1}\rm\cite[Definition 2.3]{B2012}
Let $M$ be a real symmetric matrix whose rows and columns are indexed by $X=\{1, 2,\ldots,n\}$. Let $\{X_{1},X_{2},\ldots,X_{t}\}$ be a partition of $X$. The characteristic matrix C is the $n\times m$ matrix whose $jth$ column is the characteristic vector of $X_j(j=1, 2,\ldots,m)$. Let $M$ be partitioned according to $\pi$ as

$$ M_{\pi} = \begin{bmatrix} 

    M_{11} & M_{12} & \dots & M_{1m} \\
    
    M_{21} & M_{22} & \ldots & M_{2m}   \\
    
    \vdots & \vdots & \ddots & \vdots \\

    M_{m1} & M_{m2} & \ldots & M_{mm} \\
    
\end{bmatrix}.$$
where $M_{ij}$ denotes the submatrix (block) of $M$ formed by rows in $X_i$ and the columns in $X_j$.If $q_{ij}$ denotes the average row sum of $M_{ij}$, then the matrix $Q^M=[(q^M_{ij})]_{m\times m}$ is called the quotient matrix of $M$. If the row sum of each block $M_{ij}$ is a constant, then the partition $\pi$ is called an equitable partition.
\end{definition}

\begin{definition}
Let $\mathcal{H}=\big(V(\mathcal{H}), E(\mathcal{H})\big)$ be a hypergraph and $V_{1},V_{2},\ldots,V_{t}$ be a partition of the vertex set $V(\mathcal{H})$. Then the set $\{V_{1},V_{2},\ldots,V_{t}\}$ is said to form an {\em equitable vertex partition (simply equitable partition)} if for each $r,s\in \{1,2,\ldots, t\}$ and for any $i\in V_{r}$, 
\begin{center}
$q^{\mathcal{S}}_{rs}=\sum\limits_{j;j\in V{s}}\mathcal{S}(\mathcal{H})_{ij}$,
\end{center} where $q_{rs}$ are constant depending on $r$ and $s$. 
\end{definition}
\begin{remark}
For an equitable partition, we define the Seidel quotient matrix $Q^\mathcal{S}(\mathcal{H})=[(q^{\mathcal{S}}_{rs})]_{t\times t}$.
\end{remark}

\begin{lemma}\rm\cite[Lemma 2.3.1]{B2012}\label{P1}
Let $Q^M$ be the quotient matrix of any square matrix $M$ corresponding to an equitable partition. Then the spectrum of $M$ contains the spectrum of $Q^M$.
\end{lemma}

The following theorem gives the Seidel spectrum of $\mathcal{C}^3_{m,n}$, where $m,n\geq2$. 

\begin{theorem}\label{t4}
Let $\mathcal{C}^{3}_{m,n}$ be the complete $3$-uniform bipartite hypergraph with $m,n\geq2$. Then Seidel spectrum of $\mathcal{C}^{3}_{m,n}$ is 
    
$\small{\Bigg\{\overbrace{2n-1,\ldots,2n-1}^{(m-1)},\overbrace{2m-1,\ldots,2m-1}^{(n-1)},\frac{3m+3n-4mn-2\pm\sqrt{16m^3n+n^2+2mn(8n^2-40n+49)+m^2(32n^2-80n+1)}}{2}\Bigg\}}$.
\end{theorem}
\pr Let $V(\mathcal{C}^{3}_{m,n})=\{u_1,u_2,\ldots,u_{m},u_{m+1},\ldots,u_{m+n}\}$ be the vertex set of $\mathcal{C}^{3}_{m,n}$ and $V_1=\{u_1,u_2,\ldots,u_m\}$, $V_2=\{u_{m+1},u_{m+2},\ldots,u_{m+n}\}$ be the bipartition of $\mathcal{C}^{3}_{m,n}$. Now the Seidel matrix of $\mathcal{C}^3_{m,n}$ is

$$ \mathcal{S}(\mathcal{C}^{3}_{m,n}) = \begin{bmatrix} 
(1-2n) Y_{m\times m} & \big(1-2(m+n-2)\big)X_{m\times n} \\
\big(1-2(m+n-2)\big)X^t_{m\times n} &  (1-2m) Y_{n\times n}
\end{bmatrix},$$

where $Y_{m\times m}$ and $Y_{n\times n}$ are the matrix whose diagonal entries are 0 and each non-diagonal entry is 1, and $X_{m\times n}$ is the matrix whose all entries are equal to $1$.

Let us denote $\phi_i\in\mathbb{R}^{m+n}$ be the $(m+n)$-vector with $i$-th position 1 and else 0. Let  $\mathcal{P}_{1}=\{\phi_{1,i}=\phi_1-\phi_i:i=2,3,4,\ldots,m\}$ and $\mathcal{P}_{2}=\{\phi_{m+1,j}=\phi_{m+1}-\phi_j:j=m+2,m+3,\ldots,m+n\}$.

Therefore, we have the following: 
\begin{center}
$\mathcal{S}(\mathcal{C}^{3}_{m,n})\phi_{1,i}=(2n-1)\phi_{1,i}$, $\forall\phi_{1,i}\in\mathcal{P}_1$
\end{center}

\begin{center}
$\hspace{.8cm}\mathcal{S}(\mathcal{C}^{3}_{m,n})\phi_{m+1,j}=(2m-1)\phi_{m+1,j}$, $\forall\phi_{m+1,j}\in\mathcal{P}_2$.
\end{center}

The sets of vectors $\{\phi_{1,i}:i=2,3,\ldots,m\}$ and $\{\phi_{m+1,j}:j=m+2,\ldots,m+n\}$ are linearly independent. Therefore, $2n-1$ and $2m-1$ are the Seidel eigenvalues of $\mathcal{C}^3_{m,n}$ with multiplicity at least $m-1$ and $n-1$, respectively.

We use the equitable partition to find the other Seidel eigenvalues of $\mathcal{C}^3_{m,n}$. For $\mathcal{C}^3_{m,n}$, the vertex partition, i.e., $\{V_{1},V_{2}\}$, forms an equitable vertex partition of $\mathcal{S}(\mathcal{C}^{3}_{m,n})$.

Let $Q^{\mathcal{S}}(\mathcal{C}^3_{m,n})=[q^{\mathcal{S}}_{rs}]_{2\times2}$ be the quotient matrix of $\mathcal{S}(\mathcal{C}^3_{m,n})$. Then the values of $q^{\mathcal{S}}_{11},q^{\mathcal{S}}_{12},q^{\mathcal{S}}_{21},q^{\mathcal{S}}_{22}$ will be the following:

\begin{align*}
&q^{\mathcal{S}}_{11}=\sum\limits_{i\in V_1}(1-2c_{ii})=(m-1)(1-2n)~\forall~i\in V_1;\\
&q^{\mathcal{S}}_{12}=\sum\limits_{j\in V_2}(1-2c_{ij})=n\big(1-2(m+n-2)\big)~\forall~i\in V_1;\\
&q^{\mathcal{S}}_{21}=\sum\limits_{j\in V_1}(1-2c_{ij})=m\big(1-2(m+n-2)\big)~\forall~i\in V_2;\\
&q^{\mathcal{S}}_{22}=\sum\limits_{i\in V_2}(1-2c_{ii})=(n-1)(1-2m)~\forall~i\in V_2.
\end{align*}

The quotient matrix is $Q^{\mathcal{S}}(\mathcal{C}^{3}_{m,n})=\begin{bmatrix} (m-1)(1-2n) & n\big(1-2(n+m-2)\big)\\ m\big(1-2(n+m-2)\big) & (n-1)(1-2m)\end{bmatrix}$. 

The eigenvalues of $Q^{\mathcal{S}}(\mathcal{C}^{3}_{m,n})$ are $\frac{3m+3n-4mn-2\pm\sqrt{16m^3n+n^2+2mn(8n^2-40n+49)+m^2(32n^2-80n+1)}}{2}$. Using Lemma \ref{P1}, we have that every eigenvalue of $Q^{\mathcal{S}}(\mathcal{C}^{3}_{m,n})$ is an eigenvalue of $\mathcal{S}(\mathcal{C}^{3}_{m,n})$.

Hence, the Seidel spectrum of $\mathcal{C}^{3}_{m,n}$ is 

$\small{\Bigg\{\overbrace{2n-1,\ldots,2n-1}^{(m-1)},\overbrace{2m-1,\ldots,2m-1}^{(n-1)},\frac{3m+3n-4mn-2\pm\sqrt{16m^3n+n^2+2mn(8n^2-40n+49)+m^2(32n^2-80n+1)}}{2}\Bigg\}}$.\qe\\\\
\textbf{Observation:} There are two types of hyperedge in the complete $3$-uniform bipartite hypergraph $\mathcal{C}^3_{m,n}=(V_1, V_2, E)$. Type-I: hyperedge containing one vertex from $V_1$ and two vertices from $V_2$, Type-II: hyperedge containing two vertices from $V_1$ and one vertex from $V_2$. The following example demonstrates that the Seidel spectrum of $\mathcal{C}^3_{m,n}-e_1$ and $\mathcal{C}^3_{m,n}-e_2$ may not be equal, where $e_1$ and $e_2$ are Type-I and Type-II hyperedges, respectively.

Let us consider the complete $3$-uniform bipartite hypergraphs $\mathcal{C}^3_{3,6}=(V_1,V_2, E)$. Let $V_1=\{1,2,3\}$ and $V_2=\{4,5,6,7,8,9\}$. Let $e_{1}=\{1,2,4\}$. Then the Seidel spectrum of $\mathcal{C}^3_{3,6}-e_1$ is $\{-80.291,10.252,28.040,5,5,5,5,9,13\}$. Let $e_{2}=\{1,4,5\}$. Then the Seidel spectrum $\mathcal{C}^3_{4,5}-e_2$ is $\{-80.634,27.076,13.297,10.463,5.994,2.805,5,5,11\}$. Thus, the Seidel spectra of $\mathcal{C}^3_{3,6}-e_1$ and $\mathcal{C}^3_{3,6}-e_2$ are different. So there are two types of hyperedge deletion in $\mathcal{C}^3_{m,n}$. 

\begin{remark}\label{r1} Let $\xi_1(\rho)=-\rho^5+\Big(-5+7m+5n-4mn\Big)\rho^4+\Big(46+4m-18m^2+4n-22mn+2m^2n+4m^3n-8n^2-6mn^2+4m^2n^2+4mn^3\Big)\rho^3+\Big(286-234m-46m^2+36m^3-206n+62mn-10m^2n+52m^3n-16m^4n+8n^2-10mn^2+64m^2n^2-24m^3n^2+4n^3+40mn^3-24m^2n^3-8mn^4\Big)\rho^2+\Big(83+76m-366m^2+248m^3-40m^4-228n-90mn+238m^2n+20m^3n-88m^4n+16m^5n-40n^2+246mn^2+132m^2n^2-184m^3n^2+48m^4n^2+40n^3-28mn^3-160m^2n^3+48m^3n^3-16mn^4+32m^2n^4\Big)\rho+\Big(-921+2387m-2322m^2+1012m^3-168m^4+873n-2250mn+1994m^2n-732m^3n+72m^4n+16m^5n-408n^2+1274mn^2-680m^2n^2-160m^3n^2+160m^4n^2-32m^5n^2+84n^3-368mn^3+40m^2n^3+160m^3n^3-32m^4n^3+24mn^4+32m^2n^4-32m^3n^4\Big)$ and $\rho_i$ be the roots that satisfy the equation $\xi_1(\rho)=0$, $i=1,2,3,4,5$.
\end{remark}

The following theorem gives the Seidel spectrum of $\mathcal{C}^3_{m,n}-e$ of Type-I.

\begin{theorem}\label{t7}
Let $\mathcal{C}^{3}_{m,n}=(V_1,V_2,E)$ be the complete $3$-uniform bipartite hypergraph with $m\geq2$ and $n\geq3$. Let $e\in E$ be a hyperedge containing one vertex from $V_1$ and two vertices from $V_2$. Then the Seidel spectrum of $\mathcal{C}^{3}_{m,n}-e$ is
\begin{center}
$\Bigg\{\overbrace{2n-1,2n-1\ldots,2n-1}^{(m-2)},\overbrace{2m-1,2m-1\ldots,2m-1}^{(n-3)},\rho_i(\mathcal{C}^3_{m,n}-e),i=1,2,3,4,5\Bigg\},$
\end{center}
where $\rho_i(\mathcal{C}^3_{m,n}-e) $ satisfy the equation $\xi_1(\rho)=0 $, for all $i=1,2,3,4,5 $.
\end{theorem}
\pr Let $V(\mathcal{C}^{3}_{m,n})=\{u_1,u_2,\ldots,u_{m},u_{m+1},u_{m+2},\ldots,u_{m+n}\}$ be the vertex set of $\mathcal{C}^{3}_{m,n}$ and $V_1=\{u_1,u_2,\ldots,u_m\}$, $V_2=\{u_{m+1},u_{m+2},\ldots,u_{m+n}\}$ be the bipartition of $\mathcal{C}^{3}_{m,n}$.

Since the hyperedge $e$ contains one vertex from $V_1$ and two vertices from $V_2$. Without loss of generality, let us consider $e=\{u_{1},u_{m+1},u_{m+2}\}$.

Then the Seidel matrix of $\mathcal{C}^{3}_{m,n}-e$ is the following:

\[
\resizebox{1\textwidth}{!}{
$\mathcal{S}(\mathcal{C}^3_{m,n}-e)=\begin{bmatrix} 
     0 & 1-2n & \hdots & 1-2n  & 1-2(m+n-3) & 1-2(m+n-3)& \hdots & 1-2(m+n-2)\\ 
    
     1-2n & 0 & \hdots &  1-2n & 1-2(m+n-2) & 1-2(m+n-2) & \hdots & 1-2(m+n-2) \\

    \vdots & \vdots & \ddots  &\vdots  & \vdots & \vdots & \ddots  & \vdots \\

    1-2n & 1-2n & \hdots & 0 & 1-2(m+n-2) & 1-2(m+n-2) & \hdots & 1-2(m+n-2)  \\

   1-2(m+n-3) & 1-2(m+n-2) & \hdots & 1-2(m+n-2)  & 0 & 1-2(m-1) & \hdots & 1-2m  \\

    1-2(m+n-3) & 1-2(m+n-2) & \hdots & 1-2(m+n-2) & 1-2(m-1) & 0 & \hdots & 1-2m \\

      \vdots & \vdots & \ddots  &\vdots & \vdots & \vdots & \ddots & \vdots \\

     1-2(m+n-2) & 1-2(m+n-2) & \hdots & 1-2(m+n-2) & 1-2m & 1-2m & \hdots  & 0 
\end{bmatrix}$}
\]

Let us denote $\sigma_i\in\mathbb{R}^{m+n}$ be the $(m+n)$-vector with $i$-th position 1 and else 0. Let  $\mathcal{M}_{1}=\{\sigma_{2,i}=\sigma_2-\sigma_i:i=3,4,\ldots,m\}$ and $\mathcal{M}_{2}=\{\sigma_{m+3,j}=\sigma_{m+3}-\sigma_j:j=m+4,m+5,\ldots,m+n\}$.

Therefore we have the following relationships: 
\begin{center}
$\mathcal{S}(\mathcal{C}^{3}_{m,n}-e)\sigma_{2,i}=(2n-1)\sigma_{2,i}$, $\forall\sigma_{2,i}\in\mathcal{M}_1$.
\end{center}

\begin{center}
$\hspace{.7cm}\mathcal{S}(\mathcal{C}^{3}_{m,n}-e)\sigma_{m+3,j}=(2m-1)\sigma_{m+3,j}$, $\forall\sigma_{m+3,j}\in\mathcal{M}_2$.
\end{center}

The sets of vectors $\{\sigma_{2,i}:i=3,4\ldots,m\}$ and $\{\sigma_{m+3,j}:j=m+4,m+5,\ldots,m+n\}$ are linearly independent. Therefore, $2n-1$ and $2m-1$ are the Seidel eigenvalues of $\mathcal{C}^3_{m,n}-e$ with multiplicity at least $m-2$ and $n-3$, respectively.

Now we have at least total $m-2+n-3=m+n-5$ Seidel eigenvalues of $\mathcal{C}^3_{m,n}-e$. We use equitable vertex partition to find the rest of the five Seidel eigenvalues. Let us taking $\mathcal{O}_{1}=\{u_{1}\}$, $\mathcal{O}_{2}=\{u_{m+1}\}$, $\mathcal{O}_{3}=\{u_{m+2}\}$, $\mathcal{O}_{4}=\{u_2,u_3,\ldots,u_m\}$, $\mathcal{O}_{5}=\{u_{m+3},\ldots,u_{m+n}\}$. Then the set $\{\mathcal{O}_1,\mathcal{O}_{2},\mathcal{O}_{3},\mathcal{O}_{4},\mathcal{O}_{5}\}$ forms an equitable vertex partition of $\mathcal{C}^{3}_{m,n}-e$.

Let $Q^{\mathcal{S}}(\mathcal{C}^{3}_{m,n}-e)=[(q^{\mathcal{S}}_{ij})]_{5\times5}$ be the quotient matrix of $\mathcal{S}(\mathcal{C}^3_{m,n}-e)$. Then the values of $q^\mathcal{S}_{ij}$ for $1\leq i,j\leq5$ are following:

\begin{align*}
&q^{\mathcal{S}}_{11}=q^{\mathcal{S}}_{22}=q^{\mathcal{S}}_{33}=0,~q^{\mathcal{S}}_{44}=(m-2)(1-2n)~and~q^{\mathcal{S}}_{55}=(n-3)(1-2m); \\
\vspace{.5cm}
&q^{\mathcal{S}}_{12}=q^{\mathcal{S}}_{21}=1-2(m+n-3); \\
\vspace{.5cm}
&q^{\mathcal{S}}_{13}=q^{\mathcal{S}}_{31}=1-2(m+n-3); \\
\vspace{.5cm}
&q^{\mathcal{S}}_{14}=(m-1)(1-2n)~and~q^{\mathcal{S}}_{41}=1-2n; \\
\vspace{.5cm}
&q^{\mathcal{S}}_{15}=(n-2)\big(1-2(m+n-2)\big)~ and~q^{\mathcal{S}}_{15}=1-2(m+n-2); \\
\vspace{.5cm}
&q^{\mathcal{S}}_{23}=q^{\mathcal{S}}_{32}=1-2(m+n-3); \\
\vspace{.5cm}
&q^{\mathcal{S}}_{24}=(m-1)\big(1-2(m+n-2)\big)~and~q^{\mathcal{S}}_{42}=1-2(m+n-3); \\
\vspace{.5cm}
&q^{\mathcal{S}}_{25}=(n-2)(1-2m)~and ~q^{\mathcal{S}}_{52}=1-2m ; \\
\vspace{.5cm}
&q^{\mathcal{S}}_{34}=(m-1)\big(1-2(m+n-2)\big)~and ~q^{\mathcal{S}}_{43}=1-2(m+n-3) ; \\
\vspace{.5cm}
&q^{\mathcal{S}}_{35}=(n-2)(1-2m)~and ~q^{\mathcal{S}}_{53}=1-2m ; \\
\vspace{.5cm}
&q^{\mathcal{S}}_{45}=(n-2)\big(1-2(m+n-2)\big)~and~q^{\mathcal{S}}_{54}=(m-1)\big(1-2(m+n-2)\big).
\end{align*}

Then the quotient matrix $Q^{\mathcal{S}}(\mathcal{C}^{3}_{m,n}-e) $ is defined in the Definition (\ref{d1}) as follows:

\[
\resizebox{1.02\textwidth}{!}{
$\begin{bmatrix} 
0 & 1-2(m+n-3) & 1-2(m+n-3) & (m-1)(1-2n) & (n-2)\big(1-2(m+n-2)\big) \\
    
1-2(m+n-3) & 0 & 1-2(m-1) & (m-1)\big(1-2(m+n-2)\big) & (n-2)(1-2m) \\
    
1-2(m+n-3) & 1-2(m-1) & 0 & (m-1)\big(1-2(m+n-2)\big) & (n-2)(1-2m)  \\
    
1-2n & 1-2(m+n-2) & 1-2(m+n-2) & (m-2)(1-2n) & (n-2)\big(1-2(m+n-2)\big)  \\

1-2 (m+n-2) & 1-2 m & 1-2 m & (m-1)\big(1-2(m+n-2)\big) & (n-3)(1-2 m)
\end{bmatrix}$}
\]

Using MATLAB and Remark \ref{r1}, the characteristic equation of the matrix $Q^{\mathcal{S}}(\mathcal{C}^{3}_{m,n}-e)$ is $\xi_1(\rho)=0$. By Lemma \ref{P1}, the roots of the polynomial $\xi_1(\rho) $ are in the Seidel spectrum of $\mathcal{C}^{3}_{m,n}-e$.

Hence, the Seidel spectrum of $\mathcal{C}^{3}_{m,n}-e $ is

$\Bigg\{\overbrace{2n-1,2n-1\ldots,2n-1}^{(m-2)},\overbrace{2m-1,2m-1\ldots,2m-1}^{(n-3)},\rho_i(\mathcal{C}^3_{m,n}-e),i=1,2,3,4,5\Bigg\},$

where $\rho_i(\mathcal{C}^3_{m,n}-e) $ satisfy the equation $\xi_1(\rho)=0 $, for all $i=1,2,3,4,5 $.\qe
\begin{remark}\label{r2} Let $\xi_2(\tau)=-\tau^5+\Big(-5+5m+7n-4mn\Big)\tau^4+\Big(46+4m-8m^2+4n-22mn-6m^2n+4m^3n-18n^2+2mn^2+4m^2n^2+4mn^3\Big)\tau^3+\Big(286-206m+8m^2+4m^3-234n+62mn-10m^2n+40m^3n-8m^4n-46n^2-10mn^2+64m^2n^2-24m^3n^2+36n^3+52mn^3-24m^2n^3-16mn^4\Big)\tau^2+\bigg(83-228m-40m^2+40m^3+76n-90mn+246m^2n-28m^3n-16m^4n-366n^2+238mn^2+132m^2n^2-160m^3n^2+32m^4n^2+248n^3+20mn^3-184m^2n^3+48m^3n^3-40n^4-88mn^4+48m^2n^4+16mn^5\bigg)\tau+\bigg(-921+873m-408m^2+84m^3+2387n-2250mn+1274m^2n-368m^3n+24m^4n-2322n^2+1994mn^2-680m^2n^2+40m^3n^2+32m^4n^2+1012n^3-732mn^3-160m^2n^3+160m^3n^3-32m^4n^3-168n^4+72mn^4+160m^2n^4-32m^3n^4+16mn^5-32m^2n^5\bigg)$ and $\tau_i$ be the roots that satisfy the equation $\xi_2(\tau)=0$ for all $i=1,2,3,4,5$.
\end{remark}

The following theorem gives the Seidel spectrum of $\mathcal{C}^3_{m,n}-e$ of Type-II.

\begin{theorem}\label{t9}
Let $\mathcal{C}^{3}_{m,n}=(V_1,V_2,E)$ be the complete $3$-uniform bipartite hypergraph with $m\geq3$ and $n\geq2$. Let $e\in E$ be a hyperedge containing two vertices from $V_1$ and one vertex from $V_2$. Then the Seidel spectrum of $\mathcal{C}^{3}_{m,n}-e$ is 
\begin{center}
$\Bigg\{\overbrace{2n-1,2n-1,\ldots,2n-1}^{(m-3)},\overbrace{2m-1,2m-1,\ldots,2m-1}^{(n-2)},\tau_i(\mathcal{C}^3_{m,n}-e),i=1,2,3,4,5\Bigg\},$
\end{center}
where $\tau_i(\mathcal{C}^3_{m,n}-e)$ satisfy the equation $\xi_2(\tau)=0$ for all $i=1,2,3,4,5$.
\end{theorem}
\pr Let $V(\mathcal{C}^{3}_{m,n})=\{u_1,u_2,\ldots,u_{m+n}\}$ be the vertex set of $\mathcal{C}^{3}_{m,n}$ and $V_1=\{u_1,u_2,\ldots,u_m\}$, $V_2=\{u_{m+1},u_{m+2},\ldots,u_{m+n}\}$ be the bipartition of $\mathcal{C}^{3}_{m,n}$.

Since the hyperedge $e$ contains two vertices from $V_1$ and one vertex from $V_2$. Without loss of generality, let us consider $e=\{u_{1},u_{2},u_{m+1}\}.$ Then the Seidel matrix of $\mathcal{C}^{3}_{m,n}-e$ is the following:

\[
\resizebox{1\textwidth}{!}{
$\mathcal{S}(\mathcal{C}^3_{m,n}-e)=\begin{bmatrix} 
     0 & 1-2(n-1) & \hdots & 1-2n  & 1-2(m+n-3) & 1-2(m+n-2)& \hdots & 1-2(m+n-2)\\ 
    
     1-2(n-1) & 0 & \hdots &  1-2n & 1-2(m+n-3) & 1-2(m+n-2) & \hdots & 1-2(m+n-2) \\

    \vdots & \vdots & \ddots  &\vdots  & \vdots & \vdots & \ddots  & \vdots \\

    1-2n & 1-2n & \hdots & 0 & 1-2(m+n-2) & 1-2(m+n-2) & \hdots & 1-2(m+n-2)  \\

   1-2(m+n-3) & 1-2(m+n-3) & \hdots & 1-2(m+n-2)  & 0 & 1-2(m-1) & \hdots & 1-2m  \\

    1-2(m+n-2) & 1-2(m+n-2) & \hdots & 1-2(m+n-2) & 1-2m & 0 & \hdots & 1-2m \\

      \vdots & \vdots & \ddots  &\vdots & \vdots & \vdots & \ddots & \vdots \\

     1-2(m+n-2) & 1-2(m+n-2) & \hdots & 1-2(m+n-2) & 1-2m & 1-2m & \hdots  & 0 
    \end{bmatrix}$}
\]

Let us denote $\sigma_i\in\mathbb{R}^{m+n}$ be the $(m+n)$-vector with $i$-th position 1 and else 0. Let  $\mathcal{N}_{1}=\{\sigma_{3,i}=\sigma_3-\sigma_i:i=4,5,\ldots,m\}$ and $\mathcal{N}_{2}=\{\sigma_{m+2,j}=\sigma_{m+2}-\sigma_j:j=m+3,m+4,\ldots,m+n\}$.

Therefore we have the following relationships: 
\begin{center}
    $\mathcal{S}(\mathcal{C}^{3}_{m,n}-e)\sigma_{3,i}=(2n-1)\sigma_{2,i}$, $\forall\sigma_{2,i}\in\mathcal{N}_1$.
\end{center}

\begin{center}
    $\hspace{.7cm}\mathcal{S}(\mathcal{C}^{3}_{m,n}-e)\sigma_{m+2,j}=(2m-1)\sigma_{m+2,j}$, $\forall\sigma_{m+2,j}\in\mathcal{N}_2$.
\end{center}

The sets of vectors $\{\sigma_{3,i}:i=4,5\ldots,m\}$ and $\{\sigma_{m+2,j}:j=m+3,m+4,\ldots,m+n\}$ are linearly independent. Therefore, $2n-1$ and $2m-1$ are the Seidel eigenvalues of $\mathcal{C}^3_{m,n}-e$ with multiplicity at least $m-3$ and $n-2$, respectively.

Now we have at least total $m-2+n-3=m+n-5$ Seidel eigenvalues of $\mathcal{C}^3_{m,n}-e$. We use equitable vertex partition to find the rest of the five Seidel eigenvalues. Let us taking $\mathcal{O}'_{1}=\{u_{1}\}$, $\mathcal{O}'_{2}=\{u_{2}\}$, $\mathcal{O}'_{3}=\{u_{m+1}\}$, $\mathcal{O}'_{4}=\{u_3,u_4,\ldots,u_m\}$, $\mathcal{O}'_{5}=\{u_{m+2},\ldots,u_{m+n}\}$. Then the set $\{\mathcal{O}'_1,\mathcal{O}'_{2},\mathcal{O}'_{3},\mathcal{O}'_{4},\mathcal{O}'_{5}\}$ forms an equitable vertex partition of $\mathcal{S}(\mathcal{C}^{3}_{m,n}-e)$.

Let $Q^{\mathcal{S}}(\mathcal{C}^{3}_{m,n}-e)=[(q^{\mathcal{S}}_{ij})]_{5\times5}$ be the quotient matrix of $\mathcal{S}(\mathcal{C}^3_{m,n}-e)$. Then the values of $q^\mathcal{S}_{ij}$ for $1\leq i,j\leq5$ are following:

\begin{align*}
&q^{\mathcal{S}}_{11}=q^{\mathcal{S}}_{22}=q^{\mathcal{S}}_{33}=0,~q^{\mathcal{S}}_{44}=(m-3)(1-2n)~and~q^{\mathcal{S}}_{55}=(n-2)(1-2m); \\
\vspace{1cm}
&q^{\mathcal{S}}_{12}=q^{\mathcal{S}}_{21}=1-2(n-1); \\
\vspace{1cm}
&q^{\mathcal{S}}_{13}=q^{\mathcal{S}}_{31}=1-2(m+n-3); \\
\vspace{1cm}
&q^{\mathcal{S}}_{14}=(m-2)(1-2n)~and~q^{\mathcal{S}}_{41}=1-2n; \\
\vspace{1cm}
&q^{\mathcal{S}}_{15}=(n-1)\big(1-2(m+n-2)\big)~ and~q^{\mathcal{S}}_{15}=1-2(m+n-2);\\
\vspace{1cm}
&q^{\mathcal{S}}_{23}=q^{\mathcal{S}}_{32}=1-2(m+n-3); \\
\vspace{1cm}
&q^{\mathcal{S}}_{24}=(m-2)(1-2n)~and~q^{\mathcal{S}}_{42}=1-2n; \\
\vspace{1cm}
&q^{\mathcal{S}}_{25}=(n-1)\big(1-2(m+n-2)\big)~and ~q^{\mathcal{S}}_{52}=1-2(m+n-2) ; \\
\vspace{1cm}
&q^{\mathcal{S}}_{34}=(m-2)\big(1-2(m+n-2)\big)~and ~q^{\mathcal{S}}_{43}=1-2(m+n-3) ; \\
\vspace{1cm}
&q^{\mathcal{S}}_{35}=(n-1)(1-2m)~and ~q^{\mathcal{S}}_{53}=1-2m ; \\
\vspace{1cm}
&q^{\mathcal{S}}_{45}=(n-1)\big(1-2(m+n-2)\big)~and ~q^{\mathcal{S}}_{54}=(m-1)\big(1-2(m+n-2)\big).
\end{align*}

Then the quotient matrix $Q^{\mathcal{S}}(\mathcal{C}^{3}_{m,n}-e)$ is defined in the Definition (\ref{d1}) as follows:

\[
\resizebox{1\textwidth}{!}{
$\begin{bmatrix} 
0 & 1-2(n-1) & 1-2(m+n-3) & (m-2)(1-2n) & (n-1)\big(1-2(m+n-2)\big)\\
    
1-2(n-1) & 0 & 1-2(m+n-3) & (m-2)(1-2n) & (n-1)\big(1-2(m+n-2)\big) \\
    
1-2(m+n-3) & 1-2(m+n-3) & 0 & (m-2)\big(1-2(m+n-2)\big) & (n-1)(1-2 m) \\
    
1-2n & 1-2n & 1-2(m+n-2) & (m-3)(1-2n) & (n-1)\big(1-2(m+n-2)\big) \\

1-2(m+n-2) & 1-2(m+n-2) & 1-2m & (m-2)\big(1-2(m+n-2)\big) & (n-2)(1-2m)
    \end{bmatrix}$}
\]

Using MATLAB and Remark \ref{r2}, the characteristic equation of the matrix $Q^{\mathcal{S}}(\mathcal{C}^{3}_{m,n}-e)$ is $\xi_2(\tau)=0$. By Lemma \ref{P1}, the roots of the polynomial $\xi_2(\tau) $ are in the Seidel spectrum of $\mathcal{C}^{3}_{m,n}-e$.

Hence, the Seidel spectrum of $\mathcal{C}^{3}_{m,n}-e $ is

$\Bigg\{\overbrace{2n-1,2n-1,\ldots,2n-1}^{(m-3)},\overbrace{2m-1,2m-1,\ldots,2m-1}^{(n-2)},\tau_i(\mathcal{C}^3_{m,n}-e),i=1,2,3,4,5\Bigg\},$

where $\xi_2(\tau_i)=0 $ and $\tau_i\in\mathbb{R} $, $\forall i=1,2,3,4,5 $.\qe

The following lemma says that the number of negative and positive Seidel eigenvalues of $\mathcal{C}^{3}_{m,n}$ are preserved when an arbitrary single hyperedge is deleted.

\begin{lemma}\label{TTT}
Let $\mathcal{C}^{3}_{m,n}=(V_1,V_2,E)$ be the complete $3$-uniform bipartite hypergraph with $m,n\geq3$ and let $e$ be an hyperedge of $\mathcal{C}^{3}_{m,n}$. Then $\mathcal{C}^3_{m,n}-e$ has exactly one negative Seidel eigenvalue, and the remaining Seidel eigenvalues are positive.
\end{lemma}
\pr There are two cases.

Case-I: The hyperedge $e$ contains one vertex from $V_1$ and two vertices from $V_2$. From Remark (\ref{r1}), the roots of the equation $\xi_1(\rho)=0$ are $\rho_i(\mathcal{C}^3_{m,n}-e)$, $i=1,2,3,4,5$.

Using relation between roots and coefficients of $\xi_1(\rho)$ and MATLAB, we have the following signs:
\begin{align*}
&\sum\limits_{i=1}^5\rho_i(\mathcal{C}^3_{m,n}-e)<0\hspace{.1cm}\forall m,n\geq3\\
&\sum\limits_{i\neq j}(\rho_i(\mathcal{C}^3_{m,n}-e)\rho_j(\mathcal{C}^3_{m,n}-e))<0\hspace{.1cm}\forall m,n\geq3\\
&\sum\limits_{i\neq j\neq h}\Big(\rho_i(\mathcal{C}^3_{m,n}-e)\rho_j(\mathcal{C}^3_{m,n}-e)\rho_h(\mathcal{C}^3_{m,n}-e)\Big)<0\hspace{.1cm}\forall m,n\geq3\\
&\sum\limits_{i\neq j\neq h\neq g}\Big(\rho_i(\mathcal{C}^3_{m,n}-e)\rho_j(\mathcal{C}^3_{m,n}-e)\rho_h(\mathcal{C}^3_{m,n}-e)\rho_g(\mathcal{C}^3_{m,n}-e)\Big)<0\hspace{.1cm}\forall m,n\geq3\\
&\prod\limits_{i=1}^5\rho_i(\mathcal{C}^3_{m,n}-e)<0\hspace{.1cm}\forall m,n\geq3
\end{align*}

The signs in the sequence of $\xi_1(-\rho)$ coefficients are $+-----$. There is one sign variation in $\xi_1(-\rho)$, by Descartes’ Rule of Signs the number of negative roots of $\xi_1(\rho)=0$ is exactly one. Since $\prod\limits_{i=1}^5\rho_i(\mathcal{C}^3_{m,n}-e)\neq0$, so the number of positive roots of $\xi_1(\rho)=0$ is exactly four.

Hence, using Theorem (\ref{t7}), the proof follows.

Case-II: The hyperedge $e$ contains two vertices from $V_1$ and one vertex from $V_2$. If we use Theorem (\ref{t9}), the rest of the proof is the same as Case-I.\qe

\begin{lemma}\rm\cite[Page 538, 8.4.P14]{H2013}\label{P2}
Let $A=[a_{ij}]_{n\times n}$ and $B=[b_{ij}]_{n\times n}$ be two non-negative matrices, and suppose that $A$ is irreducible. Then $\lambda(A+B)>\lambda(A)$ if $B\neq0$, where $\lambda(A+B)$ and $\lambda(A)$ are two spectral radii of $A+B,A$, respectively.
\end{lemma}

The following theorem provides an answer to \textbf{Problem 1}.
\begin{theorem}\label{E}
Let $\mathcal{C}^{3}_{m,n}=(V_1,V_2,E)$ be the complete $3$-uniform bipartite hypergraph with $m,n\geq3$. Let the hyperedge be $e\in E$. Then \begin{center}
$\mathcal{E_S}(\mathcal{C}^{3}_{m,n})>\mathcal{E_S}(\mathcal{C}^{3}_{m,n}-e)$.
\end{center}
\end{theorem}
\pr Let $e$ be the hyperedge that has been removed from $\mathcal{C}^3_{m,n}=(V_1,V_2,E)$. There are two cases.

Case-I: The hyperedge $e$ contains one vertex from $V_1$ and two vertices from $V_2$. Let $t_{i}(\mathcal{C}^3_{m,n})$ be the Seidel eigenvalues of $\mathcal{C}^3_{m,n}$, $i=1,2,3,4,5$.\\ 
Let $t_1(\mathcal{C}^3_{m,n})=\frac{3m+3n-4mn-2+\sqrt{16m^3n+n^2+2mn(8n^2-40n+49)+m^2(32n^2-80n+1)}}{2}$, $t_2(\mathcal{C}^3_{m,n})=2m-1$,\\\\ 
$t_3(\mathcal{C}^3_{m,n})=2m-1$, $t_4=2n-1$, $t_5(\mathcal{C}^3_{m,n})=\frac{3m+3n-4mn-2-\sqrt{16m^3n+n^2+2mn(8n^2-40n+49)+m^2(32n^2-80n+1)}}{2}$.

Then $t_i(\mathcal{C}^3_{m,n})>0$, $i=1,2,3,4$ and $t_5(\mathcal{C}^3_{m,n})<0$. Also, $\sum\limits_{i=1}^5t_i(\mathcal{C}^3_{m,n})=-5+7m+5n-4mn$.

Therefore, we have \begin{equation}\label{1}
\mathcal{E_S}(\mathcal{C}^3_{m,n})=(2n-1)(m-2)+(2m-1)(n-3)+t_1(\mathcal{C}^3_{m,n})+t_2(\mathcal{C}^3_{m,n})+t_3(\mathcal{C}^3_{m,n})+t_4(\mathcal{C}^3_{m,n})-t_5(\mathcal{C}^3_{m,n}).
\end{equation}

From Lemma (\ref{TTT}), there is only one negative Seidel eigenvalue of $\mathcal{C}^3_{m,n}-e$. Therefore, using Theorem (\ref{t7}), let us consider that $\rho_5(\mathcal{C}^3_{m,n}-e)<0$ and $\rho_i(\mathcal{C}^3_{m,n}-e)>0$, $i=1,2,3,4$. Also $\sum\limits_{i=1}^5\rho_i(\mathcal{C}^3_{m,n}-e)=-5+7m+5n-4mn$ \big(from Remark (\ref{r1})\big).

Therefore, we have \begin{multline}\label{2}
\mathcal{E_S}(\mathcal{C}^3_{m,n}-e)=(2n-1)(m-2)+(2m-1)(n-3)+\rho_1(\mathcal{C}^3_{m,n}-e)+\rho_2(\mathcal{C}^3_{m,n}-e)+\rho_3(\mathcal{C}^3_{m,n}-e)\\+\rho_4(\mathcal{C}^3_{m,n}-e)-\rho_5(\mathcal{C}^3_{m,n}-e).
\end{multline}

Now $-\mathcal{S}(\mathcal{C}^3_{m,n})=-\mathcal{S}(\mathcal{C}^3_{m,n}-e)+B$, where $B$ is a non-negative matrix. Since the matrix $I-\mathcal{S}(\mathcal{C}^3_{m,n}-e)$ is a positive matrix, it implies that $-\mathcal{S}(\mathcal{C}^3_{m,n}-e)$ is an irreducible matrix, where $I$ is the $(m+n)\times(m+n)$ identity matrix. Notice that $-\mathcal{S}(\mathcal{C}^3_{m,n})$ and $-\mathcal{S}(\mathcal{C}^3_{m,n}-e)$ have exactly one positive eigenvalue, which are $-t_5(\mathcal{C}^3_{m,n})$ and $-\rho_5(\mathcal{C}^3_{m,n}-e)$, respectively.

Using Lemma (\ref{P2}) $\implies -t_5(\mathcal{C}^3_{m,n})>-\rho_5(\mathcal{C}^3_{m,n}-e)$,

$\hspace{2.9cm}\implies- t_5(\mathcal{C}^3_{m,n})=-\rho_5(\mathcal{C}^3_{m,n}-e)+a$,where $a>0$

By writing in MATLAB, the polynomial in Remark (\ref{r1}) can be factored as $\xi_1(\rho)=\big(\rho-(2m-3)\big)\xi'_1(\rho)$, where $\xi'_1(\rho)$ is a $4$-degree polynomial.

Let $\rho_2(\mathcal{C}^3_{m,n}-e)=2m-3$. Then $2m-1>2m-3\implies t_2(\mathcal{C}^3_{m,n})=\rho_2(\mathcal{C}^3_{m,n}-e)+2$.

Without loss of generality, we shall compare between $t_1(\mathcal{C}^3_{m,n})$ and $\rho_1(\mathcal{C}^3_{m,n}-e)$; $t_3(\mathcal{C}^3_{m,n})$ and $\rho_3(\mathcal{C}^3_{m,n}-e)$; $t_4(\mathcal{C}^3_{m,n})$ and $\rho_4(\mathcal{C}^3_{m,n}-e)$. There are eight subcases.\\\\
Subcase(i): Let $t_1(\mathcal{C}^3_{m,n})\geq\rho_1(\mathcal{C}^3_{m,n}-e)$, $t_3(\mathcal{C}^3_{m,n})\geq\rho_3(\mathcal{C}^3_{m,n}-e)$, $t_4(\mathcal{C}^3_{m,n})\geq\rho_4(\mathcal{C}^3_{m,n}-e)$.\\\\
Then we have the following:

$t_1(\mathcal{C}^3_{m,n})=\rho_1(\mathcal{C}^3_{m,n}-e)+b$, $t_3(\mathcal{C}^3_{m,n})=\rho_3(\mathcal{C}^3_{m,n}-e)+c$, $t_4(\mathcal{C}^3_{m,n})=\rho_4(\mathcal{C}^3_{m,n}-e)+d$, where $b,c,d\geq0$.\\\\
From equation (\ref{1}) and (\ref{2}), we have,\\\\
$\mathcal{E_S}(\mathcal{C}^3_{m,n})=(m-2)(2n-1)+(n-3)(2m-1)+\rho_1(\mathcal{C}^3_{m,n}-e)+\rho_2(\mathcal{C}^3_{m,n}-e)+\rho_3(\mathcal{C}^3_{m,n}-e)+\rho_4(\mathcal{C}^3_{m,n}-e)-\rho_5(\mathcal{C}^3_{m,n}-e)+(a+b+c+d+2)$,\\\\
$\mathcal{E_S}(\mathcal{C}^3_{m,n}-e)=(m-2)(2n-1)+(n-3)(2m-1)+\rho_1(\mathcal{C}^3_{m,n}-e)+\rho_2(\mathcal{C}^3_{m,n}-e)+\rho_3(\mathcal{C}^3_{m,n}-e)+\rho_4(\mathcal{C}^3_{m,n}-e)-\rho_5(\mathcal{C}^3_{m,n}-e)$.\\\\ 
Since $b,c,d\geq 0$ and $a>0$, so $\mathcal{E_S}(\mathcal{C}^3_{m,n})>\mathcal{E_S}(\mathcal{C}^3_{m,n}-e)$. \\\\
Subcase(ii): Let $t_1(\mathcal{C}^3_{m,n})\geq\rho_1(\mathcal{C}^3_{m,n}-e)$, $t_3(\mathcal{C}^3_{m,n})\geq\rho_3(\mathcal{C}^3_{m,n}-e)$, $t_4(\mathcal{C}^3_{m,n})\leq\rho_4(\mathcal{C}^3_{m,n}-e)$.\\\\
Then we have the following:

$t_1(\mathcal{C}^3_{m,n})=\rho_1(\mathcal{C}^3_{m,n}-e)+b$, $t_3(\mathcal{C}^3_{m,n})=\rho_3(\mathcal{C}^3_{m,n}-e)+c$, $t_4(\mathcal{C}^3_{m,n})=\rho_4(\mathcal{C}^3_{m,n}-e)-d$, where $b,c,d\geq0$.\\\\
Now $t_1(\mathcal{C}^3_{m,n})+t_2(\mathcal{C}^3_{m,n})+t_3(\mathcal{C}^3_{m,n})+t_4(\mathcal{C}^3_{m,n})+t_5(\mathcal{C}^3_{m,n})=\Big(\rho_1(\mathcal{C}^3_{m,n}-e)+b\Big)+\Big(\rho_2(\mathcal{C}^3_{m,n}-e)+2\Big)+\Big(\rho_3(\mathcal{C}^3_{m,n}-e)+c\Big)+\Big(\rho_4(\mathcal{C}^3_{m,n}-e)-d\Big)+\Big(\rho_5(\mathcal{C}^3_{m,n}-e)-a\Big)$\\\\
$\implies -a+b+c-d+2=0$.\\\\
From equation (\ref{1}) and (\ref{2}), we have,\\\\
$\mathcal{E_S}(\mathcal{C}^3_{m,n})=(m-2)(2n-1)+(n-3)(2m-1)+\rho_1(\mathcal{C}^3_{m,n}-e)+\rho_2(\mathcal{C}^3_{m,n}-e)+\rho_3(\mathcal{C}^3_{m,n}-e)+\rho_4(\mathcal{C}^3_{m,n}-e)-\rho_5(\mathcal{C}^3_{m,n}-e)+(a+b+c-d+2)$,\\\\
$\mathcal{E_S}(\mathcal{C}^3_{m,n}-e)=(m-2)(2n-1)+(n-3)(2m-1)+\rho_1(\mathcal{C}^3_{m,n}-e)+\rho_2(\mathcal{C}^3_{m,n}-e)+\rho_3(\mathcal{C}^3_{m,n}-e)+\rho_4(\mathcal{C}^3_{m,n}-e)-\rho_5(\mathcal{C}^3_{m,n}-e)$.\\\\
Suppose $a+b+c-d+2\leq0$ for all $a>0$ and $b,c,d\geq0$. \\\\
Then $b+c-d+2\leq-a\implies -a+b+c-d+2\leq-2a<0$, which is a contradiction.\\\\ Therefore 
$a+b+c-d+2> 0$. Hence $\mathcal{E_S}(\mathcal{C}^3_{m,n})>\mathcal{E_S}(\mathcal{C}^3_{m,n}-e)$. \\\\
Subcase(iii): Let $t_1(\mathcal{C}^3_{m,n})\geq\rho_1(\mathcal{C}^3_{m,n}-e)$, $t_3(\mathcal{C}^3_{m,n})\leq\rho_3(\mathcal{C}^3_{m,n}-e)$, $t_4(\mathcal{C}^3_{m,n})\geq\rho_4(\mathcal{C}^3_{m,n}-e)$.\\\\
Then we have the following: 

$t_1(\mathcal{C}^3_{m,n})=\rho_1(\mathcal{C}^3_{m,n}-e)+b$,
$t_3(\mathcal{C}^3_{m,n})=\rho_3(\mathcal{C}^3_{m,n}-e)-c$,
$t_4(\mathcal{C}^3_{m,n})=\rho_4(\mathcal{C}^3_{m,n}-e)+d$, where $b,c,d\geq0$.\\\\
Now $t_1(\mathcal{C}^3_{m,n})+t_2(\mathcal{C}^3_{m,n})+t_3(\mathcal{C}^3_{m,n})+t_4(\mathcal{C}^3_{m,n})+t_5(\mathcal{C}^3_{m,n})=\Big(\rho_1(\mathcal{C}^3_{m,n}-e)+b\Big)+\Big(\rho_2(\mathcal{C}^3_{m,n}-e)+2\Big)+\Big(\rho_3(\mathcal{C}^3_{m,n}-e)-c\Big)+\big(\rho_4(\mathcal{C}^3_{m,n}-e)+d\Big)+\Big(\rho_5(\mathcal{C}^3_{m,n}-e)-a\Big)$\\\\
$\implies -a+b-c+d+2=0$.\\\\
From equation (\ref{1}) and (\ref{2}), we have,\\\\
$\mathcal{E_S}(\mathcal{C}^3_{m,n})=(m-2)(2n-1)+(n-3)(2m-1)+\rho_1(\mathcal{C}^3_{m,n}-e)+\rho_2(\mathcal{C}^3_{m,n}-e)+\rho_3(\mathcal{C}^3_{m,n}-e)+\rho_4(\mathcal{C}^3_{m,n}-e)-\rho_5(\mathcal{C}^3_{m,n}-e)+(a+b-c+d+2)$,\\\\
$\mathcal{E_S}(\mathcal{C}^3_{m,n}-e)=(m-2)(2n-1)+(n-3)(2m-1)+\rho_1(\mathcal{C}^3_{m,n}-e)+\rho_2(\mathcal{C}^3_{m,n}-e)+\rho_3(\mathcal{C}^3_{m,n}-e)+\rho_4(\mathcal{C}^3_{m,n}-e)-\rho_5(\mathcal{C}^3_{m,n}-e)$.\\\\
Suppose $a+b-c+d+2\leq0$, for all $a>0$ and $b,c,d\geq0$. \\\\
Then $b-c+d+2\leq-a\implies-a+b-c+d+2\leq-2a<0$, which is a contradiction. \\\\
Therefore 
$a+b-c+d+2> 0$. Hence $\mathcal{E_S}(\mathcal{C}^3_{m,n})>\mathcal{E_S}(\mathcal{C}^3_{m,n}-e)$.\\\\
Subcase(iv): Let $t_1(\mathcal{C}^3_{m,n})\geq\rho_1(\mathcal{C}^3_{m,n}-e)$, $t_3(\mathcal{C}^3_{m,n})\leq\rho_3(\mathcal{C}^3_{m,n}-e)$, $t_4(\mathcal{C}^3_{m,n})\leq\rho_4(\mathcal{C}^3_{m,n}-e)$.\\\\
Then we have the following:

$t_1(\mathcal{C}^3_{m,n})=\rho_1(\mathcal{C}^3_{m,n}-e)+b$,
$t_3(\mathcal{C}^3_{m,n})=\rho_3(\mathcal{C}^3_{m,n}-e)-c$,
$t_4(\mathcal{C}^3_{m,n})=\rho_4(\mathcal{C}^3_{m,n}-e)-d$, where $b,c,d\geq0$\\\\
Now $t_1(\mathcal{C}^3_{m,n})+t_2(\mathcal{C}^3_{m,n})+t_3(\mathcal{C}^3_{m,n})+t_4(\mathcal{C}^3_{m,n})+t_5(\mathcal{C}^3_{m,n})=\Big(\rho_1(\mathcal{C}^3_{m,n}-e)+b\Big)+\Big(\rho_2(\mathcal{C}^3_{m,n}-e)+2\Big)+\Big(\rho_3(\mathcal{C}^3_{m,n}-e)-c\Big)+\Big(\rho_4(\mathcal{C}^3_{m,n}-e)-d\Big)+\Big(\rho_5(\mathcal{C}^3_{m,n}-e)-a\Big)$\\\\
$\implies -a+b-c-d+2=0$.\\\\
From equation (\ref{1}) and (\ref{2}), we have,\\\\
$\mathcal{E_S}(\mathcal{C}^3_{m,n})=(m-2)(2n-1)+(n-3)(2m-1)+\rho_1(\mathcal{C}^3_{m,n}-e)+\rho_2(\mathcal{C}^3_{m,n}-e)+\rho_3(\mathcal{C}^3_{m,n}-e)+\rho_4(\mathcal{C}^3_{m,n}-e)-\rho_5(\mathcal{C}^3_{m,n}-e)+(a+b-c-d+2)$,\\\\
$\mathcal{E_S}(\mathcal{C}^3_{m,n}-e)=(m-2)(2n-1)+(n-3)(2m-1)+\rho_1(\mathcal{C}^3_{m,n}-e)+\rho_2(\mathcal{C}^3_{m,n}-e)+\rho_3(\mathcal{C}^3_{m,n}-e)+\rho_4(\mathcal{C}^3_{m,n}-e)-\rho_5(\mathcal{C}^3_{m,n}-e)$. \\\\
Suppose $a+b-c-d+2\leq0$, for all $a>0$ and $b,c,d\geq0$. \\\\
Then $b-c-d+2\leq-a$
 $\implies -a+b-c-d+2\leq-2a<0$, which is a contradiction. \\\\
Therefore 
$a+b-c-d+2> 0$. Hence $\mathcal{E_S}(\mathcal{C}^3_{m,n})>\mathcal{E_S}(\mathcal{C}^3_{m,n}-e)$.\\\\
Subcase(v): Let $t_1(\mathcal{C}^3_{m,n})\leq\rho_1(\mathcal{C}^3_{m,n}-e)$, $t_3(\mathcal{C}^3_{m,n})\geq\rho_3(\mathcal{C}^3_{m,n}-e)$, $t_4(\mathcal{C}^3_{m,n})\geq\rho_4(\mathcal{C}^3_{m,n}-e)$.\\\\
Then we have the following: 

$t_1(\mathcal{C}^3_{m,n})=\rho_1(\mathcal{C}^3_{m,n}-e)-b$, $t_3(\mathcal{C}^3_{m,n})=\rho_3(\mathcal{C}^3_{m,n}-e)+c$,
$t_4(\mathcal{C}^3_{m,n})=\rho_4(\mathcal{C}^3_{m,n}-e)+d$, where $b,c,d\geq0$.\\\\
Now $t_1(\mathcal{C}^3_{m,n})+t_2(\mathcal{C}^3_{m,n})+t_3(\mathcal{C}^3_{m,n})+t_4(\mathcal{C}^3_{m,n})+t_5(\mathcal{C}^3_{m,n})=\Big(\rho_1(\mathcal{C}^3_{m,n}-e)-b\Big)+\Big(\rho_2(\mathcal{C}^3_{m,n}-e)+2\Big)+\Big(\rho_3(\mathcal{C}^3_{m,n}-e)+c\Big)+\Big(\rho_4(\mathcal{C}^3_{m,n}-e)+d\Big)+\Big(\rho_5(\mathcal{C}^3_{m,n}-e)-a\Big)$\\\\
$\implies -a-b+c+d+2=0$.\\\\
From equation (\ref{1}) and (\ref{2}), we have,\\\\
$\mathcal{E_S}(\mathcal{C}^3_{m,n})=(m-2)(2n-1)+(n-3)(2m-1)+\rho_1(\mathcal{C}^3_{m,n}-e)+\rho_2(\mathcal{C}^3_{m,n}-e)+\rho_3(\mathcal{C}^3_{m,n}-e)+\rho_4(\mathcal{C}^3_{m,n}-e)-\rho_5(\mathcal{C}^3_{m,n}-e)+(a-b+c+d+2)$,\\\\
$\mathcal{E_S}(\mathcal{C}^3_{m,n}-e)=(m-2)(2n-1)+(n-3)(2m-1)+\rho_1(\mathcal{C}^3_{m,n}-e)+\rho_2(\mathcal{C}^3_{m,n}-e)+\rho_3(\mathcal{C}^3_{m,n}-e)+\rho_4(\mathcal{C}^3_{m,n}-e)-\rho_5(\mathcal{C}^3_{m,n}-e)$.\\\\
Suppose $a-b+c+d+2\leq0$, for all $a>0$ and $b,c,d\geq0$. \\\\
Then $\hspace{-.1cm}-b+c+d+2\leq-a\implies\hspace{-.1cm}-a-b+c+d+2\leq-2a<0$, which is a contradiction. \\\\
Therefore 
$a-b+c+d+2> 0$. Hence $\mathcal{E_S}(\mathcal{C}^3_{m,n})>\mathcal{E_S}(\mathcal{C}^3_{m,n}-e)$.\\\\
Subcase(vi): Let $t_1(\mathcal{C}^3_{m,n})\leq\rho_1(\mathcal{C}^3_{m,n}-e)$, $t_3(\mathcal{C}^3_{m,n})\geq\rho_3(\mathcal{C}^3_{m,n}-e)$ and $t_4(\mathcal{C}^3_{m,n})\leq\rho_4(\mathcal{C}^3_{m,n}-e)$.\\\\
Then we have the following: 

$t_1(\mathcal{C}^3_{m,n})=\rho_1(\mathcal{C}^3_{m,n}-e)-b$,
$t_3(\mathcal{C}^3_{m,n})=\rho_3(\mathcal{C}^3_{m,n}-e)+c$,
$t_4(\mathcal{C}^3_{m,n})=\rho_4(\mathcal{C}^3_{m,n}-e)-d$, where $b,c,d\geq0$.\\\\
Now $t_1(\mathcal{C}^3_{m,n})+t_2(\mathcal{C}^3_{m,n})+t_3(\mathcal{C}^3_{m,n})+t_4(\mathcal{C}^3_{m,n})+t_5(\mathcal{C}^3_{m,n})=\Big(\rho_1(\mathcal{C}^3_{m,n}-e)-b\Big)+\Big(\rho_2(\mathcal{C}^3_{m,n}-e)+2\Big)+\Big(\rho_3(\mathcal{C}^3_{m,n}-e)+c\Big)+\Big(\rho_4(\mathcal{C}^3_{m,n}-e)-d\Big)+\Big(\rho_5(\mathcal{C}^3_{m,n}-e)-a\Big)$\\\\
$\implies -a-b+c-d+2=0$.\\\\
From equation (\ref{1}) and (\ref{2}), we have,\\\\
$\mathcal{E_S}(\mathcal{C}^3_{m,n})=(m-2)(2n-1)+(n-3)(2m-1)+\rho_1(\mathcal{C}^3_{m,n}-e)+\rho_2(\mathcal{C}^3_{m,n}-e)+\rho_3(\mathcal{C}^3_{m,n}-e)+\rho_4(\mathcal{C}^3_{m,n}-e)-\rho_5(\mathcal{C}^3_{m,n}-e)+(a-b+c-d+2)$,\\\\
$\mathcal{E_S}(\mathcal{C}^3_{m,n}-e)=(m-2)(2n-1)+(n-3)(2m-1)+\rho_1(\mathcal{C}^3_{m,n}-e)+\rho_2(\mathcal{C}^3_{m,n}-e)+\rho_3(\mathcal{C}^3_{m,n}-e)+\rho_4(\mathcal{C}^3_{m,n}-e)-\rho_5(\mathcal{C}^3_{m,n}-e)$.\\\\
Suppose $a-b+c-d+2\leq0$, for all $a>0$ and $b,c,d\geq0$. \\\\
Then $\hspace{-.1cm}-b+c-d+2\leq-a\implies\hspace{-.1cm}-a-b+c-d+2\leq-2a<0$, which is a contradiction. \\\\
Therefore 
$a-b+c-d+2> 0$. Hence $\mathcal{E_S}(\mathcal{C}^3_{m,n})>\mathcal{E_S}(\mathcal{C}^3_{m,n}-e)$.\\\\
Subcase(vii): Let $t_1(\mathcal{C}^3_{m,n})\leq\rho_1(\mathcal{C}^3_{m,n}-e)$, $t_3(\mathcal{C}^3_{m,n})\leq\rho_3(\mathcal{C}^3_{m,n}-e)$ and $t_4(\mathcal{C}^3_{m,n})\geq\rho_4(\mathcal{C}^3_{m,n}-e)$.\\\\
Then we have the following:

$t_1(\mathcal{C}^3_{m,n})=\rho_1(\mathcal{C}^3_{m,n}-e)-b$,
$t_3(\mathcal{C}^3_{m,n})=\rho_3(\mathcal{C}^3_{m,n}-e)-c$,
$t_4(\mathcal{C}^3_{m,n})=\rho_4(\mathcal{C}^3_{m,n}-e)+d$, where $b,c,d\geq0$.\\\\
Now $t_1(\mathcal{C}^3_{m,n})+t_2(\mathcal{C}^3_{m,n})+t_3(\mathcal{C}^3_{m,n})+t_4(\mathcal{C}^3_{m,n})+t_5(\mathcal{C}^3_{m,n})=\Big(\rho_1(\mathcal{C}^3_{m,n}-e)-b\Big)+\Big(\rho_2(\mathcal{C}^3_{m,n}-e)+2\Big)+\Big(\rho_3(\mathcal{C}^3_{m,n}-e)-c\Big)+\Big(\rho_4(\mathcal{C}^3_{m,n}-e)+d\Big)+\Big(\rho_5(\mathcal{C}^3_{m,n}-e)-a\Big)$\\\\
$\implies-a-b-c+d+2=0$.\\\\
From equation (\ref{1}) and (\ref{2}), we have,\\\\
$\mathcal{E_S}(\mathcal{C}^3_{m,n})=(m-2)(2n-1)+(n-3)(2m-1)+\rho_1(\mathcal{C}^3_{m,n}-e)+\rho_2(\mathcal{C}^3_{m,n}-e)+\rho_3(\mathcal{C}^3_{m,n}-e)+\rho_4(\mathcal{C}^3_{m,n}-e)-\rho_5(\mathcal{C}^3_{m,n}-e)+(a-b-c+d+2)$,\\\\
$\mathcal{E_S}(\mathcal{C}^3_{m,n}-e)=(m-2)(2n-1)+(n-3)(2m-1)+\rho_1(\mathcal{C}^3_{m,n}-e)+\rho_2(\mathcal{C}^3_{m,n}-e)+\rho_3(\mathcal{C}^3_{m,n}-e)+\rho_4(\mathcal{C}^3_{m,n}-e)-\rho_5(\mathcal{C}^3_{m,n}-e)$.\\\\
Suppose $a-b-c+d+2\leq0$, for all $a>0$ and $b,c,d\geq0$. \\\\
Then $\hspace{-.1cm}-b-c+d+2\leq-a\implies\hspace{-.1cm}-a-b-c+d+2\leq-2a<0$, which is a contradiction. \\\\
Therefore 
$a-b-c+d+2> 0$. Hence $\mathcal{E_S}(\mathcal{C}^3_{m,n})>\mathcal{E_S}(\mathcal{C}^3_{m,n}-e)$.\\\\
Subcase(viii): Let $t_1(\mathcal{C}^3_{m,n})\leq\rho_1(\mathcal{C}^3_{m,n}-e)$, $t_3(\mathcal{C}^3_{m,n})\leq\rho_3(\mathcal{C}^3_{m,n}-e)$ and $t_4(\mathcal{C}^3_{m,n})\leq\rho_4(\mathcal{C}^3_{m,n}-e)$.\\\\
Then we have the following:

$t_1(\mathcal{C}^3_{m,n})=\rho_1(\mathcal{C}^3_{m,n}-e)-b$,
$t_3(\mathcal{C}^3_{m,n})=\rho_3(\mathcal{C}^3_{m,n}-e)-c$,
$t_4(\mathcal{C}^3_{m,n})=\rho_4(\mathcal{C}^3_{m,n}-e)-d$, where $b,c,d\geq0$.\\\\
Now $t_1(\mathcal{C}^3_{m,n})+t_2(\mathcal{C}^3_{m,n})+t_3(\mathcal{C}^3_{m,n})+t_4(\mathcal{C}^3_{m,n})+t_5(\mathcal{C}^3_{m,n})=\Big(\rho_1(\mathcal{C}^3_{m,n}-e)-b\Big)+\Big(\rho_2(\mathcal{C}^3_{m,n}-e)+2\Big)+\Big(\rho_3(\mathcal{C}^3_{m,n}-e)-c\Big)+\Big(\rho_4(\mathcal{C}^3_{m,n}-e)-d\Big)+\Big(\rho_5(\mathcal{C}^3_{m,n}-e)-a\Big)$\\\\
$\implies-a-b-c-d+2=0$.\\\\
From equation (\ref{1}) and (\ref{2}), we have,\\\\
$\mathcal{E_S}(\mathcal{C}^3_{m,n})=(m-2)(2n-1)+(n-3)(2m-1)+\rho_1(\mathcal{C}^3_{m,n}-e)+\rho_2(\mathcal{C}^3_{m,n}-e)+\rho_3(\mathcal{C}^3_{m,n}-e)+\rho_4(\mathcal{C}^3_{m,n}-e)-\rho_5(\mathcal{C}^3_{m,n}-e)+(a-b-c-d+2)$,\\\\
$\mathcal{E_S}(\mathcal{C}^3_{m,n}-e)=(m-2)(2n-1)+(n-3)(2m-1)+\rho_1(\mathcal{C}^3_{m,n}-e)+\rho_2(\mathcal{C}^3_{m,n}-e)+\rho_3(\mathcal{C}^3_{m,n}-e)+\rho_4(\mathcal{C}^3_{m,n}-e)-\rho_5(\mathcal{C}^3_{m,n}-e)$. \\\\
Suppose $a-b-c-d+2\leq0$, for all $a>0$ and $b,c,d\geq0$. \\\\
Then $\hspace{-.1cm}-b-c-d+2\leq-a\implies\hspace{-.1cm}-a-b-c-d+2\leq-2a<0$, which is a contradiction. \\\\
Therefore 
$a-b-c-d+2> 0$. Hence $\mathcal{E_S}(\mathcal{C}^3_{m,n})>\mathcal{E_S}(\mathcal{C}^3_{m,n}-e)$.

Hence, from all eight sub-cases, we finally prove that \begin{center}
$\mathcal{E_S}(\mathcal{C}^3_{m,n})>\mathcal{E_S}(\mathcal{C}^3_{m,n}-e)$.
\end{center}

Case-II: The hyperedge $e$ contains two vertices from $V_1$ and one vertex from $V_2$. The rest of the proof is the same as Case-I.\qe

\begin{definition}
The Turán hypergraph, denoted by $T(n,k,r)$, is a $k$-uniform complete multi-partite hypergraph of vertices $n$; it is formed by partitioning a set of vertices into $r$ subsets $V_1,V_2,\ldots,V_r$ with cardinality $n_1,n_2,\ldots,n_r$ such that $|n_i-n_j|\leq 1$, and then connecting $k$ vertices by a hyperedge if and only if they belong to different subsets.
\end{definition}

The following result is an immediate corollary of Theorem \ref{E}.
\begin{corollary}
For any hyperedge e of Turán hypergraph $T(n,3,2)$ with $n\geq6$, \begin{center}
$\mathcal{E_S}(T(n,3,2))>\mathcal{E_S}(T(n,3,2)-e)$.
\end{center}
\end{corollary}
\pr This is a particular case of Theorem \ref{E}.\qe

We conclude this section with an open problem:

\begin{question}
For $k\ge4$, how does the Seidel energy of the complete $k$-uniform bipartite hypergraph change under hyperedge deletion?
\end{question}

\section{Monotonicity of the Seidel energy under vertex deletion}

\begin{definition}\rm\cite{vol2009}
Let $\mathcal{H}=\big(V(\mathcal{H}),~E(\mathcal{H})\big)$ be a hypergraph and $v \in V(\mathcal{H})$ a vertex. 

\begin{itemize}
    \item The \emph{strong deletion} of $v$ removes the vertex together with all hyperedges incident to it, i.e., $$\mathcal{H} - v = \bigl(V(\mathcal{H}) \setminus \{v\},\, E'\bigr), ~~~E' = \{\, e \in E(\mathcal{H}): v \notin e \,\}.$$

    \item The \emph{weak deletion} of $v$ removes the vertex from the vertex set but retains each hyperedge after deleting $v$ from it, i.e., $$\mathcal{H}- v = 
\bigl(V(\mathcal{H}) \setminus \{v\},\, E\bigr), ~~~E' = \{\, e \setminus \{v\} : e \in E(\mathcal{H}) \,\}.$$
\end{itemize}
In the graph, strong and weak vertex deletion are the same. However, in a hypergraph, they are not the same.
\end{definition}

In 2023, the authors proved that the Seidel energy of a hypergraph decreases under weak vertex deletion.

\begin{lemma}\label{weak vertex}\rm\cite[Theorem 4.4]{LJK2023}
Let $G^{\ast}=\big(V(G^{\ast}), E(G^{\ast})\big)$ be a hypergraph of order $n$, and let $v \in V(G^{\ast})$ be any arbitrary vertex. 
Then
\[
\mathcal{E_S}(G^{\ast}) \;\ge\; \mathcal{E_S}(G^{\ast}-v).
\]
\end{lemma}

The following example demonstrates that the Seidel energy of a hypergraph does not necessarily decrease under strong vertex deletion.

\begin{example}
Let  $\mathcal{H}^*=\big(V(\mathcal{H}^*), E(\mathcal{H}^*)\big)$ be a $3$-uniform hypergraph with $V(\mathcal{H}^*)=\{1,2,3,4,5,6\}$ and $E(\mathcal{H}^*)=\Big\{\{1,2,3\},\{1, 2, 5\},\{1,3,6\},\{1,5,6\},\{2,4,6\},\{3,4,5\}\Big\}$. Let $v=4$.

Then the Seidel spectrum of $\mathcal{H}^*$ and $\mathcal{H}^*-v$ are $\Big\{-7.808,~-0.635,~1,~1,~1,~5.443\Big\}$ and $\Big\{-7.168,~-1,~-1,~0.574,~3,~5.594\Big\} $.

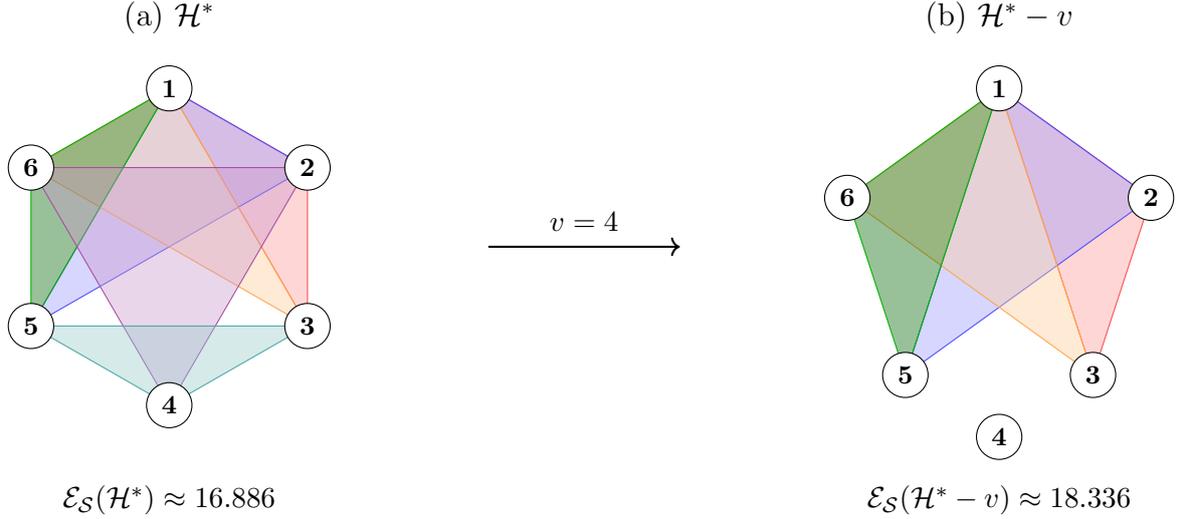
\begin{figure}[h!]
\centering
\begin{tikzpicture}[scale=1.05,
    vertex/.style={circle, draw=black, fill=white, minimum size=17pt, inner sep=0pt, font=\small\bfseries}
]

\begin{scope}[xshift=-5.2cm]
\coordinate (v1) at (90:2);
\coordinate (v2) at (30:2);
\coordinate (v3) at (-30:2);
\coordinate (v4) at (-90:2);
\coordinate (v5) at (-150:2);
\coordinate (v6) at (150:2);

\begin{scope}[fill opacity=0.4]
  \filldraw[fill=red!40, draw=red!70, rounded corners=8pt] (v1)--(v2)--(v3)--cycle;
  \filldraw[fill=blue!40, draw=blue!70, rounded corners=8pt] (v1)--(v2)--(v5)--cycle;
  \filldraw[fill=orange!40, draw=orange!70, rounded corners=8pt] (v1)--(v3)--(v6)--cycle;
  \filldraw[fill=green!40!black, draw=green!70!black, rounded corners=8pt] (v1)--(v5)--(v6)--cycle;
  \filldraw[fill=violet!40, draw=violet!70, rounded corners=8pt] (v2)--(v4)--(v6)--cycle;
  \filldraw[fill=teal!40, draw=teal!70, rounded corners=8pt] (v3)--(v4)--(v5)--cycle;
\end{scope}

\node[vertex] at (v1) {1};
\node[vertex] at (v2) {2};
\node[vertex] at (v3) {3};
\node[vertex] at (v4) {4};
\node[vertex] at (v5) {5};
\node[vertex] at (v6) {6};

\node at (0,2.9) {\large (a) $\mathcal{H}^*$};
\node at (0,-3.2) {$\mathcal{E_S}(\mathcal{H}^*)\approx 16.886$};
\end{scope}

\begin{scope}[xshift=5.2cm]
\coordinate (w1) at (90:2);
\coordinate (w2) at (18:2);
\coordinate (w3) at (-54:2);
\coordinate (w5) at (-126:2);
\coordinate (w6) at (162:2);
\coordinate (w4) at (270:2.4);

\begin{scope}[fill opacity=0.4]
  \filldraw[fill=red!40, draw=red!70, rounded corners=8pt] (w1)--(w2)--(w3)--cycle;
  \filldraw[fill=blue!40, draw=blue!70, rounded corners=8pt] (w1)--(w2)--(w5)--cycle;
  \filldraw[fill=orange!40, draw=orange!70, rounded corners=8pt] (w1)--(w3)--(w6)--cycle;
  \filldraw[fill=green!40!black, draw=green!70!black, rounded corners=8pt] (w1)--(w5)--(w6)--cycle;
\end{scope}

\node[vertex] at (w1) {1};
\node[vertex] at (w2) {2};
\node[vertex] at (w3) {3};
\node[vertex] at (w5) {5};
\node[vertex] at (w6) {6};
\node[vertex] at (w4) {4};

\node at (0,2.9) {\large (b) $\mathcal{H}^*-v$};
\node at (0,-3.2) {$\mathcal{E_S}(\mathcal{H}^*-v) \approx 18.336$};
\end{scope}

\draw[->, thick] (-1.2,0) -- (1.2,0) node[midway, above=2pt] {$v=4$};
\end{tikzpicture}
\caption{Comparison of Seidel energy of a hypergraph due to a strong vertex deletion}
\label{fig:H-hyperedge-blob-horizontal}
\end{figure}

Hence, the strong deletion of the vertex $v=4$ increases the Seidel energy, i.e. $\mathcal{E_S}(\mathcal{H}^*)<\mathcal{E_S}(\mathcal{H}^*-v)$.
\end{example}

The following theorem shows that for fixed $m\ge2$, the Seidel energy of the complete $3$-uniform bipartite hypergraph $\mathcal{C}^3_{m,n}$ decreases strictly with $n$.

\begin{theorem}\label{thm:mono}
Let $\mathcal{C}^{3}_{m,n}$ be the complete $3$-uniform bipartite hypergraph with $m,n\geq3$. Then for every fixed $m\ge2$ and all $n\ge2$,
\[
\mathcal{E_S}\!\bigl(\mathcal{C}^3_{m,n}\bigr)
\;<\;
\mathcal{E_S}\!\bigl(\mathcal{C}^3_{m,n+1}\bigr).
\]
\end{theorem}
\pr By Theorem~\ref{t4}, the Seidel spectrum of $\mathcal{C}^3_{m,n}$ consists of
\[
\overbrace{2n-1,\ldots,2n-1}^{m-1\text{ times}},\qquad
\overbrace{2m-1,\ldots,2m-1}^{n-1\text{ times}},\qquad
\rho_+(\mathcal{C}^3_{m,n}),~~\rho_-(\mathcal{C}^3_{m,n}),
\]

where $\rho_\pm(\mathcal{C}^3_{m,n})$ are the two eigenvalues of the following quotient matrix
\[
Q^{\mathcal{S}}(\mathcal{C}^3_{m,n})=
\begin{bmatrix}
(m-1)(1-2n) & n\!\left(1-2(m+n-2)\right)\\[2pt]
m\!\left(1-2(m+n-2)\right) & (n-1)(1-2m)
\end{bmatrix}.
\]

Let $T=\mathrm{tr}\,\Big(Q^{\mathcal{S}}(\mathcal{C}^3_{m,n})\big)=3m+3n-4mn-2$
and
\[
R=\det\big(Q^{\mathcal{S}}(\mathcal{C}^3_{m,n})\Big)=(m-1)(n-1)(1-2n)(1-2m)-mn\!\left(1-2(m+n-2)\right)^{\!2}\!<0\quad(m,n\ge2).
\]

Let $\rho_+(\mathcal{C}^3_{m,n})>0>\rho_-(\mathcal{C}^3_{m,n})$. Then
\[
|\rho_+(\mathcal{C}^3_{m,n})|+|\rho_-(\mathcal{C}^3_{m,n})|=\rho_+(\mathcal{C}^3_{m,n})-\rho_-(\mathcal{C}^3_{m,n})=\sqrt{T^2-4R}.
\]

Therefore the Seidel energy of $\mathcal{C}^3_{m,n}$ is
\begin{equation}\label{eq:ES-formula}
\mathcal{E_S}\!\bigl(\mathcal{C}^3_{m,n}\bigr)=(m-1)(2n-1)+(n-1)(2m-1)+\sqrt{U(m,n)}\;,
\end{equation}

where  $U(m,n)=T^2-4R=16 m^3 n+n^2+2mn(8n^2-40n+49)+m^2(32n^2-80n+1)$.

For fixed $m$, let $L(m,n):=(m-1)(2n-1)+(n-1)(2m-1)$,

and hence
\[
L(m,n+1)-L(m,n)=2(m-1)+(2m-1)=4m-3>0\qquad(m\ge2).
\]

Since the square-root function is a positive real-valued function, it is increasing. Therefore it suffices to show
\[
\Delta U(m,n):=U(m,n+1)-U(m,n)>0.
\]

Now we have
\begin{equation}\label{eq:DeltaU}
\Delta U(m,n)=16m^3+64m^2n-48m^2+48mn^2-112mn+34m+2n+1.
\end{equation}

For $n\ge2$, we have $64m^2n\ge128m^2$, $48mn^2\ge192m$, and $2n\ge4$ and substituting these into the equation~\eqref{eq:DeltaU} we get,
\[
\Delta U(m,n)
\ge16m^3+(128-48)m^2+(192-112+34)m+(4+1)
=16m^3+80m^2+114m+5>0,
\]

From equation~\eqref{eq:ES-formula} we finally get,
\[
\mathcal{E_S}\!\bigl(\mathcal{C}^3_{m,n+1}\bigr)
-\mathcal{E_S}\!\bigl(\mathcal{C}^3_{m,n}\bigr)
=\underbrace{(4m-3)}_{>0}
+\underbrace{\bigl(\sqrt{U(m,n+1)}-\sqrt{U(m,n)}\bigr)}_{>0}>0.
\]

Therefore $\mathcal{E_S}(\mathcal{C}^3_{m,n})$ is strictly increasing in~$n$.\\\\
This completes the proof.\qe 

\begin{remark}
By symmetry between $m$ and $n$ in $\mathcal{C}^3_{m,n}$,
the same argument we can prove that $\mathcal{E_S}(\mathcal{C}^3_{m,n})$ is also strictly increasing in $m$ for fixed $n$.
\end{remark}

Directly from Theorem \ref{thm:mono}, deduce that the Seidel energy of the complete 3-uniform bipartite hypergraph $\mathcal{C}^3_{m,n}$ decreases under vertex deletion for all $m,n\geq 3$. 

\begin{corollary}
Let $\mathcal{C}^{3}_{m,n}$ be the complete $3$-uniform bipartite hypergraph with $m,n\geq3$ and let $v$ be a vertex of $\mathcal{C}^{3}_{m,n}$. Then for strong vertex deletion of $v$,  
\[
\mathcal{E_S}\!\bigl(\mathcal{C}^3_{m,n}\bigr)
\;>\;
\mathcal{E_S}\!\bigl(\mathcal{C}^3_{m,n}-v\bigr).
\]
\end{corollary}

\section*{Declaration of competing interest}
There are no conflicts of interest.

\section*{Acknowledgements}
The author thanks anonymous referees for their careful review and constructive suggestions. Special thanks are due to Dr.~Swarup Kumar Panda for his insightful comments. The author expresses sincere gratitude to the Department of Mathematics (BIU) and gratefully acknowledges the postdoctoral financial support provided by Bar-Ilan University (BIU), Israel.

\section*{Data availability}
No data was used for the research described in the article.

\end{document}